\documentclass[amssymb,11pt]{amsart}
\usepackage{amscd,amssymb,euscript,amsthm}
\theoremstyle{plain}
\newtheorem{thm}{Theorem}[section] 

\newtheorem{prop}[thm]{Proposition}
\newtheorem{lem}[thm]{Lemma}
\theoremstyle{definition}
\newtheorem{defn}[thm]{Definition}
\theoremstyle{remark}
\newtheorem{rem}[thm]{Remark}

\numberwithin{equation}{section}
\newcommand{\diag}{\operatorname{diag}}

\newcommand{\id}{\operatorname{id}}
\newcommand{\ran}{\operatorname{ran}}

\newcommand{\ball}{\operatorname{ball}}
\newcommand{\dig}{\operatorname{diag\,}}
\newcommand{\aut}{\operatorname{aut}}
\newcommand{\rev}{\operatorname{Rev}}
\def\<{\left<}
\def\>{\right>}
\newcommand{\comment}[1]{}
\begin{document}
\title[Asymptotic Liftings]{the asymptotic lift of a \\ Completely Positive Map}
\author{William Arveson}
\thanks{*supported by 
NSF grant DMS-0100487} 
\address{Department of Mathematics,
University of California, Berkeley, CA 94720}
\email{arveson@math.berkeley.edu}
\maketitle

\begin{abstract} Starting with a 
unit-preserving normal completely positive  
map $L: M\to M$ acting on a von Neumann algebra - or more generally 
a dual operator system - we show that there 
is a unique reversible system $\alpha: N\to N$ (i.e., a complete 
order automorphism $\alpha$ of a dual operator system $N$) 
that captures all of the asymptotic behavior of $L$,  
called the {\em asymptotic lift} of $L$.  
This provides a noncommutative 
generalization of the Frobenius theorems  
that describe the asymptotic behavior of the sequence of powers of a stochastic 
$n\times n$ matrix.  
In cases where $M$ is a von Neumann algebra, 
the asymptotic lift is shown to be a $W^*$-dynamical system $(N,\mathbb Z)$, 
and we identify $(N,\mathbb Z)$ as the tail 
flow of the minimal dilation of $L$.  We are also able to identify the Poisson boundary 
of $L$ as the fixed algebra $N^\alpha$.

In general, we show the action of the asymptotic lift  is trivial iff 
$L$ is {\em slowly oscillating} in the sense that 
$$
\lim_{n\to\infty}\|\rho\circ L^{n+1}-\rho\circ L^n\|=0,\qquad \rho\in M_* .   
$$ 
Hence $\alpha$ is often a nontrivial automorphism of $N$.   The asymptotic lift 
of a variety of examples is calculated.  
\end{abstract}
\maketitle

\section{introduction}\label{S:it}

Throughout this paper we use the term UCP map to denote 
a normal unit-preserving completely positive 
map $L:M_1\to M_2$ of one dual operator system into another.  
While we are primarily concerned with the dynamical properties of UCP maps 
$L:M\to M$ that act on a von Neumann algebra $M$, it is 
appropriate to broaden that category to include UCP 
self-maps of more general dual operator systems.  

Stochastic $n\times n$ matrices $P = (p_{ij})$ 
describe the transition probabilities of 
$n$-state Markov chains.  The asymptotic properties of 
the sequence of powers of the transition matrix  govern the 
long-term statistical behavior of the process after initial transient  
fluctuations have died out (\cite{doobBk}, pp. 170--185).  
A stochastic $n\times n$ matrix $P=(p_{ij})$ gives rise to a 
UCP map of the commutative von Neumann algebra $\mathbb C^n$ by way of 
$$
(Px)_i=\sum_{j=1}^n p_{ij}x_j,\qquad 1\leq i\leq n, \quad 
x=(x_1,\dots,x_n)\in\mathbb C^n,  
$$
and the classical Perron-Frobenius theory provides an effective description 
of the asymptotic behavior of the sequence   
$P,P^2, P^3,\dots$ (see Section \ref{S:ex}).  

Recently, several emerging 
areas of mathematics have opened a vista of potential applications for noncommutative  
generalizations of the Frobenius theorems.  The most obvious examples are 
quantum probability, 
noncommutative dynamics \cite{arvMono}, and quantum computing \cite{kup1}.  
Such considerations led us to initiate a study  
of ``almost periodic" UCP maps on von Neumann algebras 
in \cite{arvStableI}.  Broadly speaking, those results allowed us  
to relate the asymptotic behavior of the 
powers of {\em certain} UCP maps on von Neumann algebras to the asymptotic 
behavior of 
$*$-automorphisms of other naturally associated von Neumann algebras.

However, the hypothesis of almost periodicity is very restrictive.  For example, it 
forces the spectrum of the UCP map to have a discrete component consisting of 
eigenvalues on the unit circle.  While this discrete spectrum 
appeared to be essential for the results of \cite{arvStableI}, it is 
missing in many important examples.  Nevertheless, in {\em all} of the  examples 
that we were able to penetrate, there was a ``hidden" 
W$^*$-dynamical system that shared the same asymptotic behavior.  
That led us to suspect that a different formulation 
might be possible, in which the almost periodic hypothesis 
is eliminated entirely.  This turned out to 
be true, but the new results require a complete reformulation 
in terms of UCP 
maps on dual operator systems.  

The main results of this paper 
are Theorems \ref{inThm1}, \ref{vnThm1}, \ref{ntThm1}, \ref{isThm1}.  
While these general results are 
not as sharp as those of \cite{arvStableI} when restricted 
to the case of 
almost periodic maps, we believe that 
is compensated for by the simplicity and full generality of 
the new setting.  
There are natural variations of all of the above results for one-parameter 
semigroups of normal completely positive maps that will 
be taken up elsewhere.  We also note that 
a very recent paper of St\o rmer \cite{storPos} 
complements the results of \cite{arvStableI}.

Finally, I want to thank Masaki Izumi for useful discussions about  
the properties of noncommutative Poisson boundaries during a pleasant visit to 
Kyoto in 2005.  We will make use of the noncommutative Poisson 
boundary in the proof of Theorem \ref{vnThm1} below;  and we relate it to 
the fixed algebra of the asymptotic lift in Section \ref{S:pb}.

\section{Reversible Lifts of UCP maps}\label{S:in}
Recall that an {\em operator system} is a norm-closed self-adjoint linear 
subspace $M$ of the algebra $\mathcal B(H)$ of all bounded operators on 
a Hilbert space $H$, such that the identity operator $\mathbf 1$ belongs 
to $M$.  A {\em dual} operator system is an operator system that is 
closed in the weak$^*$-topology of $\mathcal B(H)$.  We write 
$M_*$ for the predual of such an operator system $M$, namely the 
norm-closed linear subspace of the dual of $M$ consisting 
of all restrictions to $M$ of normal linear functionals on $\mathcal B(H)$; and 
since $M$ can be naturally identified with the dual of $M_*$, we refer to  
the $M_*$-topology as the weak$^*$-topology of $M$.  A {\em normal} linear 
map of dual operator systems is a linear map $L:M\to N$ that is continuous 
relative to the 
respective weak$^*$-topologies.  
Of course, such a map is the adjoint of a unique bounded linear map
of preduals, namely $\rho\in N_*\mapsto \rho\circ L\in M_*$.

While there is an abstract characterization of operator systems 
(\cite{ChE1}, Theorem 4.4), we shall not 
require the details of it here.  Only rarely do we require 
a realization $M\subseteq \mathcal B(H)$ of $M$ as a concrete 
dual operator system; but when we do, 
we require that the realization have the 
above properties - namely that $M$ is a weak$^*$-closed concrete operator 
system in $\mathcal B(H)$ and $M_*$ consists of all restrictions 
of normal linear functionals of $\mathcal B(H)$.  
For the most part, our conventions and basic terminology 
follow those of the monographs \cite{EffRu} and \cite{paulsenBk2}.  

We fix attention on the category whose objects are UCP self-maps 
$$
L:M\to M
$$ 
acting on dual operator systems $M$, and whose morphisms are 
equivariant UCP maps.  Thus,  a homomorphism from $L_1: M_1\to M_1$ to 
$L_2: M_2\to M_2$ is a UCP map 
$$
E:M_1\to M_2
$$ 
satisfying $E\circ L_1=L_2\circ E$.  
In this case we say that $L_1$ is 
a {\em lifting} of $L_2$ through $E$, or simply a {\em lift} of $L_2$. 
By an {\em automorphism} 
of a dual operator system $N$ we mean a 
UCP map $\alpha: N\to N$ having a UCP inverse $\alpha^{-1}: N\to N$.  
\begin{defn}\label{inDef0}
A {\em reversible lift} of a UCP map $L:M\to M$ is a triple $(N,\alpha,E)$ consisting 
of an automorphism $\alpha: N\to N$ of another dual operator system $N$ and a UCP map 
$E: N\to M$ satisfying $E\circ\alpha=L\circ E$.  
\end{defn}

A UCP map $L:M\to M$ has many reversible lifts, the simplest being the trivial lift 
$(\mathbb C,{\rm id}, \iota)$, where $\iota: \mathbb C\to M$ 
is the inclusion $\iota(\lambda)=\lambda\cdot\mathbf 1_M$.  
We begin by pointing out that all reversible lifts must satisfy 
a system of asymptotic inequalities.  In the usual way, $L:M\to M$ gives rise 
to a hierarchy of UCP maps $L_n: M^{(n)}\to M^{(n)}$, $n=1,2,\dots$, in which 
$M^{(n)}=M_n(\mathbb C)\otimes M$ is the $n\times n$ matrix system 
over $M$ and $L_n={\rm id}\otimes L: M^{(n)}\to M^{(n)}$ is the naturally induced 
UCP map.  Similarly, 
a reversible lifting $(N,\alpha, E)$ of  $L:M\to M$ gives rise 
to a  hierarchy of reversible liftings 
$(N^{(n)},\alpha_n, E_n)$ of 
$L_n:M^{(n)}\to M^{(n)}$, one for every $n=1,2,\dots$.  

\begin{prop}\label{inProp1}
Let $(N,\alpha,E)$ be a reversible lift of a given UCP map 
$L:M\to M$.  
For every bounded linear functional $\rho$ on $M$, the sequence of norms $\|\rho\circ L^k\|$ 
decreases with increasing $k=1,2,\dots$, and we have  
\begin{equation}\label{inEq1}
\|\rho\circ E\|\leq \lim_{k\to\infty}\|\rho\circ L^k\|.  
\end{equation}

Moreover, the inequalities (\ref{inEq1}) persist throughout the 
hierarchy of liftings of $L_n: M^{(n)}\to M^{(n)}$, as $\rho$ ranges over 
the dual of $M^{(n)}$, $n\geq 1$.  
\end{prop}

\begin{proof}Every UCP map is a contraction, hence 
$\|\rho\circ L^{k+1}\|\leq \|\rho\circ L^k\|$ for every $k\geq 0$.  
Moreover, for fixed $x\in N$ and $k\geq 0$ we can write 
$$
E(x)=E(\alpha^k(\alpha^{-k}(x)))=L^k(E(\alpha^{-k}(x))).  
$$
Since $\|E\circ\alpha^{-k}\|\leq 1$, we conclude 
that for every $\rho\in M^\prime$, 
$$
\|\rho\circ E\|=\sup_{\|x\|=1}|\rho(E(x))|=\sup_{\|x\|=1}|\rho\circ L^k(E\circ\alpha^{-k}(x))|
\leq \|\rho\circ L^k\|,
$$ 
and (\ref{inEq1}) follows 
after passing to the limit as $k\to\infty$.  
The same argument applies throughout the hierarchy of liftings 
$(N^{(n)}, \id_n\otimes\alpha, \id_n\otimes E)$, after one  
replaces $\rho$ with a bounded linear functional on $M^{(n)}$, $n=1,2,\dots$.  
\end{proof}

\begin{rem}[Nondegeneracy] \label{inRem1}  Let $(N,\alpha,E)$ be a reversible lift of 
a UCP map $L: M\to M$.  
Since $E: N\to M$ is normal, it has a pre-adjoint $E_*: M_*\to N_*$, 
defined by $E_*(\rho)=\rho\circ E$, $\rho\in M_*$.  Consider the range 
$$
E_*(M_*)=\{\rho\circ E: \rho\in M_*\}
$$ 
of this map.  
$E_*(M_*)$ is a linear subspace of $N_*$, and note that it is invariant 
under the invertible isometry $\alpha_*\in\mathcal B(N_*)$:  
\begin{equation*}
E_*(M_*)\circ\alpha\subseteq E_*(M_*).  
\end{equation*}
Indeed, that is immediate from equivariance and normality of $L$, 
since for $\rho\in M_*$ we 
have $\rho\circ E\circ\alpha=(\rho\circ L)\circ E\in M_*\circ E$.  It follows 
that the sequence $E_*(M_*)\circ \alpha^n$, $n\in\mathbb Z$,  defines a 
doubly infinite tower of subspaces 
of $N_*$
\begin{equation*}
\cdots \subseteq E_*(M_*)\circ\alpha\subseteq 
E_*(M_*)\subseteq E_*(M_*)\circ \alpha^{-1}\subseteq 
E_*(M_*)\circ\alpha^{-2}\subseteq\cdots .
\end{equation*}
A straightforward application of the Hahn-Banach theorem shows that 
this tower is norm-dense in $N_*$ if and only if for every $y\in N$ one has 
\begin{equation}\label{inEq0}
 E(\alpha^{-n}(y))=0, \ n=0,1,2,\dots \ \implies y=0.  
\end{equation}
A reversible lifting $(N,\alpha,E)$ of $L$ is said to be {\em nondegenerate} 
if condition (\ref{inEq0}) is satisfied.  

It is significant that when (\ref{inEq0}) fails, 
one can always replace $(N,\alpha,E)$ 
with a {\em nondegenerate}  
reversible lifting $(\tilde N, \tilde\alpha,\tilde E)$ 
using the following device.  We may assume that $M\subseteq \mathcal B(H)$ 
is realized as a concrete dual operator system.  Let $\tilde H=\ell^2(\mathbb Z)\otimes H$  
be the Hilbert space of all square-summable bilateral sequences from $H$ 
and define a map $\theta: N\to \mathcal B(\tilde H)$ by 
\begin{equation}\label{inEq0.1}
(\theta(y)\xi)(n)=E(\alpha^{-n}(y))\xi(n),\qquad \xi\in \tilde H, \quad n\in\mathbb Z.  
\end{equation}
$\tilde N$ is defined as the weak$^*$-closure of $\theta(N)$.  The unitary 
shift defined on $\tilde H$ by 
$U\xi(n)=\xi(n-1)$, $n\in\mathbb Z$, implements a $*$-automorphism $\tilde\alpha(X)=UXU^*$ 
of $\mathcal B(\tilde H)$ such that $\tilde\alpha(\tilde N)=\tilde N$; indeed, $\theta(N)$ is 
stable under shifts to the left or right because $\alpha$ is an automorphism of $N$.  
Note too that:
$$
(U\theta(y)U^{-1}\xi)(n)=E(\alpha^{-n+1}(y)\xi(n)=L(E(\alpha^{-n}(y)))\xi(n), 
\ \ \xi\in \tilde H, \ n\in\mathbb Z.  
$$
It follows that the map $\tilde E:\mathcal B(\tilde H)\to \mathcal B(H)$ that compresses an operator matrix in $\mathcal B(\tilde H)$ to
its $00$th component restricts to a UCP map $\tilde E:\tilde N\to M$ satisfying  
$\tilde E\circ\tilde \alpha=L\circ\tilde E$.  Thus, $(\tilde N, \tilde \alpha, \tilde E)$ 
is a reversible lifting of $L$.  It is is a homomorphic image of $(N,\alpha,E)$ in the sense 
that the UCP map $\theta:N\to \tilde N$ satisfies 
$\tilde\alpha\circ \theta=\theta\circ\alpha$ and $\tilde E\circ \theta=E$ (see Section \ref{S:hi}
for a discussion of the category of reversible liftings of $L$).  Finally, 
by examining components in the obvious way, 
one verifies directly that $(\tilde N, \tilde \alpha,\tilde E)$ is nondegenerate.  
\end{rem}

\section{Asymptotic Lifts of UCP maps}\label{S:al}

In  this section we show by a direct construction that there is a reversible lifting with favorable 
asymptotic properties and that, after degeneracies 
have been eliminated, it is {\em unique} up to natural isomorphism.

\begin{defn}\label{inDef1}
Let $L:M\to M$ be a UCP map on a dual operator system.  
An {\em asymptotic lift} of $L$ 
is a reversible lifting $(N,\alpha, E)$ of $L$ 
that satisfies nondegeneracy (\ref{inEq0}), 
such that the inequalities (\ref{inEq1}) become equalities for {\em normal} 
linear functionals throughout 
the entire matrix hierarchy 
\begin{equation}\label{inEqA}
\|\rho\circ (\id_n\otimes E)\|=\lim_{k\to\infty}\|\rho\circ (\id_n\otimes L)^k\|,\quad 
\rho\in M^{(n)}_*, \ n=1,2,\dots.    
\end{equation}
\end{defn}

We come now to a basic result.  

\begin{thm}\label{inThm1}
Every UCP map $L:M\to M$ of a dual operator system has an asymptotic lifting.  
If $(N_1, \alpha_1, E_1)$ and 
$(N_2, \alpha_2, E_2)$ are two asymptotic liftings 
for $L$, then there is a 
unique isomorphism of dual operator systems $\theta: N_1\to N_2$ such that 
$\theta\circ\alpha_1=\alpha_2\circ\theta$ and $E_2\circ\theta=E_1$.  
\end{thm}

The existence assertion of Theorem \ref{inThm1} is proved by a direct 
construction involving inverse sequences,
which are defined as follows:    

\begin{defn}\label{inDef2}
Let $L:M\to M$ be a UCP map on a dual operator system.  
By an {\em inverse sequence} for $L$ we mean a bilateral sequence $(x_n)_{n\in\mathbb Z}$ 
of elements of $M$ satisfying $\sup_n\|x_n\|<\infty$, and 
\begin{equation}\label{isEq1}
x_n=L(x_{n+1}), \qquad n\in\mathbb Z.  
\end{equation}
The set of all inverse sequences for $L$ is denoted $\mathcal S_L$, or more simply $\mathcal S$ 
when there is no cause for confusion.  
\end{defn}

\begin{rem}[Properties of Inverse sequences]\label{isRem1}
The set $\mathcal S$ of all inverse sequences for $L:M\to M$ 
is a vector space that is closed 
under pointwise involution $(x_n)\mapsto (x_n^*)$, and it contains 
all ``constant" scalar  sequences of the form 
$(\cdots,\lambda\cdot\mathbf 1, \lambda\cdot\mathbf 1,\lambda\cdot\mathbf 1,\cdots )$, 
$\lambda\in\mathbb C$.  More 
generally, the constant sequences 
$(\cdots,a,a,a,\cdots )\in \mathcal S$ 
correspond bijectively with the space of fixed elements $\{a\in M: L(a)=a\}$.  
Notice too that $\mathcal S$ is stable under shifting to the right or left; if $(x_n)_{n\in\mathbb Z}$ 
belongs to $\mathcal S$ then so does $(x_{n+k})_{n\in\mathbb Z}$ for every $k=0,\pm1, \pm2,\dots$.   

Every element $x_k$ of an inverse sequence $(x_n)$ determines all of its predecessors 
uniquely, since $x_{k-1}=L(x_k)$ and, more generally, 
$x_r=L^{k-r}(x_k)$ for all $r\leq k$.  
On the other hand, $x_k$ does {\em not} determine $x_{k+1}$ uniquely, since there 
can be many solutions $z$ of the equation $L(z)=x_k$.  If we fix a particular solution 
$z$ and replace $x_{k+1}$ with $z$ in the $k+1$\,st spot,  
then it may not be possible to solve the equation 
$z=L(w)$ for $w\in M$; and {\em in that case there is no inverse sequence 
whose $k+1$st term is the replaced element $z$ and whose $k$th term is $x_k$}.  

More generally, given an element $a\in M$, the question of whether or not 
there is an inverse sequence $(x_n)$ satisfying $x_0=a$ can be subtle.  
\end{rem}

\begin{proof}[Proof of Theorem \ref{inThm1}: Existence]
In order to construct an asymptotic lifting for $L$, 
we realize $M\subseteq \mathcal B(H)$ as a concrete weak$^*$-closed operator system.   
Let $\ell^2\otimes H$ be the Hilbert space of all sequences $n\in\mathbb Z\mapsto \xi_n\in H$ 
satisfying $\sum_n\|\xi_n\|^2<\infty$, with its usual inner product.  

We can realize the space  $\mathcal S$ of all inverse sequences for $L$ as  
an operator subspace  $N\subseteq \mathcal B(\ell^2\otimes H)$ by identifying 
an inverse sequence $(x_n)\in\mathcal S$ with the diagonal operator $D=\dig (x_n)$ defined by 
$$
(D\xi)_n=x_n\xi_n,\qquad n\in\mathbb Z.  
$$
The operator norm of $\dig (x_n)$ is given by 
$
\|\dig (x_n)\|=\sup_{n\geq 0}\|x_n\|.  
$
Since $L$ is a normal map, the relations defined by (\ref{isEq1}) are weak$^*$-closed, 
and therefore the space $N=\{\diag(x_n): (x_n)\in \mathcal S\}$ 
is closed in the weak$^*$-topology of 
$\mathcal B(\ell^2\otimes H)$.  Remark \ref{isRem1} implies that $N$ 
is self-adjoint and  
contains the identity operator of $\mathcal B(\ell^2\otimes H)$, hence  
$N$ acquires the structure of a dual operator system.

Consider the right-shift $\alpha: N\to N$ of diagonal operators 
$$
\alpha(\dig (x_n))=\dig (x_{n-1}).  
$$
Letting $U$ be the unitary bilateral shift acting on $\ell^2\otimes H$ by 
$$
(U\xi)_n=\xi_{n-1},\qquad n\in\mathbb Z, \quad \xi\in \ell^2\otimes H, 
$$
one finds that  the associated $*$-automorphism $X\mapsto UXU^*$ of 
$B(\ell^2\otimes H)$ implements the action of $\alpha$ on $N$.   
Hence $\alpha$ is a complete automorphism of 
the concrete operator system $N$.  
The natural inclusion of $H$ in $\ell^2(\mathbb Z)\otimes H$,  in which 
a vector $\xi\in H$ is identified with $\delta_0\otimes\xi\in \ell^2(\mathbb Z)\otimes H$,  
gives rise to a normal map $X\mapsto P_HX\restriction_H$ 
that restricts to a map $E: N\to M$ satisfying 
$$
E(\dig (x_n))=x_0, \qquad (x_n)\in\mathcal S.     
$$
One has 
$$
E\circ\alpha(\dig (x_n))=E(\dig (x_{n-1}))=x_{-1}=L(x_0)=L\circ E(\dig (x_n)),   
$$
so that $(N,\alpha,E)$ becomes a reversible lift of $L$.  

It remains to verify (\ref{inEq0}) and (\ref{inEqA}).  For (\ref{inEq0}), 
suppose that $X=\dig(x_k)$ satisfies $E(\alpha^{-n}(X))=x_n=0$ for $n\geq 0$.  
Since the sequence $(x_k)$ belongs to $S$, this implies that $x_j=L^{|j|}(x_0)=0$ 
for negative $j$ as well, hence $(x_k)$ is the zero sequence.  

For (\ref{inEqA}), it is enough to verify the system of 
nontrivial inequalities 
\begin{equation}\label{inEq4}
\|\rho\circ (\id_n\otimes E)\|\geq\lim_{k\to\infty}\|\rho\circ (\id_n\otimes L)^k\|,\quad 
\rho\in M^{(n)}_*, \ n=1,2,\dots.    
\end{equation}
Consider first the case $n=1$, and choose $\rho\in M_*$. 
For each $k=1,2,\dots$ we can find an element $u_k\in M$ satisfying 
$\|u_k\|=1$ and $|\rho(L^k(u_k))|=\|\rho\circ L^k\|$.  Consider the triangular array 
$s^k=(s^k_0,\dots,s^k_k)$, $k=1,2,\dots$, of elements of $M$ defined by 
$$
(s^k_0,\dots,s^k_k)=(L^k(u_k), L^{k-1}(u_k),\dots, L(u_k),u_k), \qquad k=1,2,\dots.  
$$
Each component of every one of these sequences belongs to $\ball M$, 
and $s^k_0=x_0$ for every $k\geq 1$.   Moreover, 
the $j$\,th component $s^k_j$ of any one of them is obtained from the $j+1$\,st by 
applying $L$, 
\begin{equation}\label{inEq7}
s^k_j=L(s^k_{j+1}), \qquad j=0,1,2,\dots, k-1.  
\end{equation}
Since the infinite cartesian product $\ball M\times\ball M\times\cdots$ is compact 
in its weak$^*$ product topology, 
there is a subnet $s^{k^\prime}$ 
of the sequence $s^k$ with the property that each of its components 
(note that each component is well-defined for sufficiently large $k^\prime$) 
converges weak$^*$ to an element of $\ball M$.  Hence we 
can define a single infinite sequence 
$x_0, x_1, x_2,\dots\in\ball M$ by 
$$
x_j=\lim_{k^\prime} s^{k^\prime}_j,\qquad j=0, 1,2,\dots.  
$$
Since $L$ is a normal map, the relations (\ref{inEq7}) imply that 
$x_{j+1}=L(x_j)$, $j=0,1,2,\dots$.  If we continue the sequence 
into negative integers by setting $x_k=L^{|k|}(x_0)$ for $k\leq -1$, the
result is a sequence $(x_n)\in \mathcal S$ satisfying $\sup_n\|x_n\|\leq 1$  
and $E(\dig (x_n))=x_0$.  Thus we conclude that 
$$
\|\rho\circ E\|\geq |\rho(x_0)|=\lim_{k^\prime}|\rho(L^{k^\prime}(u_{k^\prime}))|=
\lim_{k^\prime}\|\rho\circ L^{k^\prime}\|=\lim_{k\to\infty}\|\rho\circ L^k\|, 
$$   
and (\ref{inEq4}) follows.  

Notice that this argument can be repeated {\em verbatim} to establish  
(\ref{inEq4}) throughout the 
matrix hierarchy for $n\geq 2$, since 
the inverse sequences for 
$\id_n\otimes L$ are bilateral sequences $(\tilde x_k)$ whose components $\tilde x_k$ are 
$n\times n$ matrices in $M^{(n)}$ that satisfy $(\id_n\otimes L)(\tilde x_{k+1})=\tilde x_k$, 
$k\in\mathbb Z$.  
We conclude that $(N,\alpha,E)$ is an asymptotic lift of $L$.
\end{proof}

Turning now to the uniqueness issue for asymptotic lifts, we require 
the dual formulation of (\ref{inEqA}) - Lemma \ref{inLem1} 
below - 
the proof of which makes use of the following elementary result.  Since we lack 
an appropriate reference, 
we sketch a proof of the latter for completeness.  

\begin{lem}\label{inLem0}
Let $X$ be a Banach space, let $K_1\supseteq K_2\supseteq \cdots$ be a 
decreasing sequence of nonempty weak$^*$-compact convex subsets of the dual $X^\prime$
with intersection $K_\infty$.  Then for every weak$^*$-continuous linear functional 
$\rho\in X^{\prime\prime}$,
\begin{equation}\label{inEq6}
\sup\{|\rho(x)|:x\in K_n\}\downarrow\sup\{|\rho(x)|: x\in K_\infty\}, \qquad {\rm as\ } 
n\to\infty.  
\end{equation}
\end{lem}

\begin{proof}
Fix $\rho$.  The sequence of nonnegative numbers $\sup\{|\rho(x)|: x\in K_n\}$ 
obviously decreases with $n$, and its limit $\ell$ satisfies 
$$
\ell\geq\sup\{|\rho(x)|: x\in K_\infty\}.
$$   
To prove the opposite inequality choose, for every $n=1,2,\dots$, an element 
$x_n\in K_n$ such that $|\rho(x_n)|=\sup\{|\rho(x)|:x\in K_n\}$.  By compactness of 
the unit ball of $X^\prime$, there is a subnet $\{x_{n^\prime}\}$ of $\{x_n\}$ 
that converges weak$^*$ to $x_\infty$.   The limit point 
$x_\infty$ must belong to $K_\infty=\cap_nK_n$ 
because the $K_n$ decrease with $n$, and since the numbers $|\rho(x_n)|=\sup\{|\rho(x)|:x\in K_n\}$ 
converge to $\ell$, it follows from weak$^*$-continuity of $\rho$ that 
$$
|\rho(x_\infty)|=\lim_{n^\prime}|\rho(x_{n^\prime})|=\lim_{n\to\infty}|\rho(x_n)|=\ell.  
$$
Since $\ell=|\rho(x_\infty)|\leq\sup\{|\rho(x)|: x\in K_\infty\}$, the proof is complete. 
\end{proof}

\begin{lem}\label{inLem1}
Let  $L:M\to M$ be a UCP map on a dual operator system and let 
$(N,\alpha, E)$ be a reversible lift of $L$.  
The following are equivalent:
\begin{enumerate}
\item[(i)] 
Every $\rho\in M_*$ satisfies (\ref{inEqA}) 
$$
\|\rho\circ E\|=\lim_{n\to\infty}\|\rho\circ L^n\|.  
$$
\item[(ii)]
Writing $\ball X$ for the closed unit ball of a normed space $X$, we have 
\begin{equation}\label{inEq3}
E(\ball N)=\bigcap_{n=0}^\infty L^n(\ball M).  
\end{equation}
\end{enumerate}
\end{lem}

\begin{proof} (i)$\implies$(ii): Choose an element $y\in \ball N$.  
Then for every $n=0,1,2,\dots$ we can write 
$$
E(y)=E(\alpha^n\circ\alpha^{-n}(y))=L^nE(\alpha^{-n}(y))\in L^n(\ball M)
$$
since $E(\alpha^{-n}(y))\in \ball M$.  Hence $E(\ball N)\subseteq \cap_nL^n(\ball M))$.  
For the opposite inclusion, note that both sides of 
$E(\ball N)\subseteq\cap_n L^n(\ball M)$ are circled weak$^*$-compact 
convex subsets of $M$, so by a standard separation theorem it suffices 
to show that for every $\rho\in M_*$ one has 
\begin{equation}\label{inEq8}
\sup\{|\rho(x)|: x\in E(\ball N)\}=\sup\{|\rho(x)|: x\in \cap_nL^n(\ball M)\}.   
\end{equation}
The left side of (\ref{inEq8}) is $\|\rho\circ E\|$,  while  by Lemma \ref{inLem0}, 
the right side is 
$$
\lim_{n\to\infty}\{|\rho(x)|: x\in L^n(\ball M)\}=\lim_{n\to\infty}\|\rho\circ L^n\|.  
$$
An application of the hypothesis (i) now gives (\ref{inEq8}).  

(ii)$\implies$(i):  Choose $\rho\in N_*$.  For $n=1,2,\dots$ let $K_n=L^n(\ball M)$, 
and set  
$K_\infty=\cap_n K_n$.  
Lemma \ref{inLem0} implies that $\|\rho\circ L^n\|=\sup\{|\rho(y)|: y\in K_n\}$ 
decreases to $\sup\{|\rho(y)|: y\in K_\infty\}$ as $n\uparrow \infty$, while by 
(\ref{inEq3}), 
$$
\|\rho\circ E\|=\sup\{|\rho(y)|: y\in E(\ball N)\}=\sup\{|\rho(y)|: y\in K_\infty\}.     
$$ 
Thus $\lim_n\|\rho\circ L^n\|$ is identified with $\|\rho\circ E\|$, and (i) follows.  
\end{proof}

\begin{lem}\label{inLem2} Let  $L:M\to M$ be a UCP map of a dual operator system and let 
$(N,\alpha, E)$ be a reversible lift of $L$ that satisfies 
\begin{equation}\label{inEq5.1}
\|\rho\circ E\|=\lim_{n\to\infty}\|\rho\circ L^n\|,\qquad \rho\in M_*.  
\end{equation}
Let $K=\{z\in N: E(\alpha^n(z))=0, \ k\in\mathbb Z\}$.  Then for every $y\in N$ we have 
\begin{equation}\label{inEq5}
\sup_{k\in\mathbb Z}\|E(\alpha^k(y))\|=\inf_{z\in K}\|y+z\|.     
\end{equation}
\end{lem}

\begin{proof}  The inequality $\leq$ is apparent, since for  
$z\in K$ and $k\in\mathbb Z$ we have 
$$
\|E(\alpha^{k}(y))\|=\|E(\alpha^{k}(y+z))\|\leq \|y+z\|.    
$$

In order to prove $\geq$, we may assume that $\sup_k\|E(\alpha^k(y)\|=1$.  
For fixed $n=1,2,\dots$, note first that $E(\alpha^{-n}(y))\in L^p(\ball M)$ for 
every $p=1,2,\dots$.  Indeed, we have 
$
E(\alpha^{-n}(y))=L^p(E(\alpha^{-n-p}(y))
$, and $\|E(\alpha^{-n-p}(y))\|\leq 1$.  Lemma \ref{inLem1} implies that 
$\cap_p L^p(\ball M)=E(\ball N)$; and since $\alpha^n$ is an invertible isometry 
of $N$, we can 
find an element $y_n\in N$ satisfying $\|y_n\|\leq 1$ and $E(\alpha^{-n}(y_n))=E(\alpha^{-n}(y))$.  

Note that $y$ and $y_n$ satisfy the following relations 
\begin{equation}\label{inEq5.2}
E(\alpha^{-k}(y_n))=E(\alpha^{-k}(y)),\qquad k\in\mathbb Z,\quad k\leq n.  
\end{equation}
Indeed, an application of $L^{n-k}$ to both sides of $E(\alpha^{-n}(y_n))=E(\alpha^{-n}(y))$ 
leads to 
$$
E(\alpha^{-k}(y_n))=L^{n-k}E(\alpha^{-n}(y_n))=L^{n-k}E(\alpha^{-n}(y))=E(\alpha^{-k}(y)).  
$$
By weak$^*$-compactness of the unit ball of $N$, there is a subnet 
$y_{n_\alpha}$ of the sequence $y_n$ that 
converges weak$^*$ to an element $y_\infty\in  N$ such that $\|y_\infty\|\leq 1$.  Since 
each map $z\mapsto E(\alpha^{-k}(z))$ is weak$^*$-continuous on the 
unit ball of $N$, the equations (\ref{inEq5.2}) imply that 
$
E(\alpha^{-k}(y_\infty))=E(\alpha^{-k}(y)),      
$
$k\in \mathbb Z$.  
Hence $y_\infty-y$ belongs to $K$.  It follows that 
$$
\sup_{k\in\mathbb Z}\|E(\alpha^k(y))\|=1\geq \|y_\infty\| =\|y+(y_\infty-y)\|\geq \inf_{z\in K}\|y+z\|
$$
proving the inequality $\geq$ of (\ref{inEq5}).  
\end{proof}

\begin{proof}[Proof of Theorem \ref{inThm1}: Uniqueness]
Let $(N, \alpha,E)$ be the asymptotic lift of $L$ constructed 
in the above existence proof and let $(\tilde N, \tilde \alpha, \tilde E)$ be another one.    Define 
a map $\theta: y\in \tilde N\to \theta(y)\in N$ as follows: 
$$
\theta(y)=\dig (x_n), \qquad {\rm where\ } x_n=\tilde E(\tilde\alpha^{-n}(y)),\quad n\in\mathbb Z.  
$$ 
Obviously, $\theta$ is a UCP map of $\tilde N$ into $N$
satisfying 
$\theta\circ\tilde\alpha=\alpha\circ\theta$.   Lemma \ref{inLem2} implies 
that $\theta$ is an injective isometry.  Indeed, applying Lemma \ref{inLem2} repeatedly 
to the sequence of maps $\id_n\otimes\theta: \tilde N^{(n)}\to N^{(n)}$, $n=1,2,\dots$,  we find that 
$\theta$ is 
a complete isometry.  We claim that $\theta(\tilde N)=N$.  
Since $\theta$ is a weak$^*$-continuous isometry, its range is a weak$^*$-closed 
subspace of $N$, so to prove $\theta(\tilde N)=N$ it suffices to show that $\theta(\tilde N)$ 
is weak$^*$-dense in $N$.  For that, let $(x_k)$ be an inverse sequence satisfying 
$\sup_k\|x_k\|=1$.  Then 
$x_n$ belongs to $E(\ball N)$ for every $n\geq 1$; and 
by Lemma \ref{inLem1}, $E(\ball N)=\tilde E(\ball \tilde N)$.  Hence there is 
an element $y_n\in\ball\tilde N$ such that $x_n=\tilde E(\tilde \alpha^{-n}(y_n))$.  
As in the proof of (\ref{inEq5.2}), this implies $x_k=\tilde E(\tilde \alpha^{-k}(y_n))$ 
for all $k\leq n$.  These equations imply that 
$\theta(y_1),\theta(y_2),\dots$ converges component-by-component to 
$\diag (x_k)$ in the weak$^*$-topology.  
Therefore $\theta(y_n)\to\diag(x_k)$ (weak$^*$) as $n\to\infty$.  

Hence $\theta$ is an equivariant isomorphism of $\tilde N$ on $N$.  Since the zeroth component of 
$\theta(y)$ is $\tilde E(y)$, 
we also have $E\circ\theta=\tilde E$.  

Finally, we claim that 
the two requirements $\theta\circ\tilde\alpha=\alpha\circ\theta$ and $E\circ\theta=\tilde E$ 
serve to determine such a UCP map $\theta$ uniquely.  
Indeed, if 
$\theta_1, \theta_2: \tilde N\to N$ are two equivariant UCP maps satisfying 
$E\circ\theta_1=E\circ\theta_2=\tilde E$, then for every $y\in \tilde N$, 
$n\in\mathbb Z$ and 
$k=1,2$, we have 
$$
E(\alpha^n(\theta_k(y)))=E(\theta_k(\tilde\alpha^n(y)))=\tilde E(\tilde\alpha^n(y)))
$$
so that $E(\alpha^n(\theta_1(y)-\theta_2(y)))=0$, $n\in\mathbb Z$.  
Since $(N,\alpha,E)$ is nondegenerate, this implies $\theta_1(y)-\theta_2(y)=0$.  
\end{proof}

\section{The Hierarchy of Reversible lifts of $L$}\label{S:hi}

We have emphasized that the asymptotic lift of a UCP map $L:M\to M$ is characterized 
by a family of asymptotic formulas (\ref{inEqA}).  We now show that it is 
possible to characterize the asymptotic lift of $L$ in a way that 
makes no explicit reference to the asymptotic behavior of the sequence $L, L^2, L^3,\dots$, 
but rather makes use of a natural 
ordering of the set of all (equivalence classes of) nondegenerate reversible liftings of $L$.  

Throughout this brief section, $L:M\to M$ will be a fixed UCP map acting on a dual 
operator system $M$.  
We consider the category $\rev _L$ whose objects are {\em nondegenerate} reversible liftings 
$(N, \alpha,E)$ of $L$; a homomorphism from $(N_1, \alpha_1,E_1)$ to 
$(N_2, \alpha_2, E_2)$ is by definition a UCP map $\theta: N_1\to N_2$ that 
satisfies $\theta\circ\alpha_1=\alpha_2\circ\theta$ and $E_2\circ\theta=E_1$.  
Significantly, a homomorphism in this category gives rise to an embedding of 
operator systems:
\begin{prop}
A homomorphism $\theta: (N_1,\alpha_1,E_1)\to (N_2,\alpha_2,E_2)$ 
of nondegenerate reversible lifts of $L$ defines an injective 
map of $N_1$ to $N_2$.  
\end{prop}

\begin{proof}
If $y\in N_1$ satisfies $\theta(y)=0$, then for every $n\in \mathbb Z$ we 
have 
$$
E_1(\alpha_1^n(y))=E_2(\theta(\alpha_1^n(y)))=E_2(\alpha_2^n(\theta(y)))=0, 
$$
hence $y=0$ by nondegeneracy (\ref{inEq0}) of $(N_1,\alpha_1,E_1)$.  
\end{proof}

When a homomorphism $\theta:(N_1, \alpha_1,E_1)\to (N_2,\alpha_2,E_2)$ exists, we write 
$$
(N_1,\alpha_1,E_1)\leq (N_2, \alpha_2,E_2).
$$  
There is an obvious notion of isomorphism in this category, namely that there 
should exist a map $\theta$ as above which has a UCP inverse $\theta^{-1}: N_2\to N_1$.  
Obviously, both relations $\geq$ and $\leq$ hold between isomorphic elements of $\rev_L$.  
Conversely,

\begin{prop}\label{hiProp1}
Any two nondegenerate reversible liftings $(N_k,\alpha_k,E_k)$ that satisfy 
$
(N_1,\alpha_1,E_1)\leq (N_2,\alpha_2,E_2)\leq (N_1,\alpha_1,E_1)
$
are isomorphic.  
\end{prop}

\begin{proof}  By hypothesis, there are equivariant UCP maps $\theta: N_1\to  N_2$ and 
$\phi: N_2\to N_1$ such that $E_2\circ\theta=E_1$ and $E_1\circ\phi=E_2$.  Consider 
the composition $\phi\circ\theta: N_1\to N_1$.   We claim that 
$\phi\circ\theta$ is the identity map of $N_1$.  Indeed, for every 
$y\in N_1$ and $n\in\mathbb Z$, one can write 
\begin{align*}
E_1(\alpha_1^n((\phi\circ\theta)(y)))&=E_1(\alpha_1^n\circ\phi\circ\theta(y))=
E_1(\phi\circ\alpha_2^n\circ\theta(y))\\
&=E_2(\alpha_2^n\circ\theta(y))=E_2(\theta\circ\alpha_1^n(y))=E_1(\alpha_1^n(y)).  
\end{align*}
It follows that $E_1(\alpha_1^n(\phi\circ\theta(y)-y))=0$ for all $n\in\mathbb Z$.   
By the nondegeneracy hypothesis (\ref{inEq0}) we conclude that $\phi\circ\theta(y)=y$.  
By symmetry,  $\theta\circ\phi$ is the identity map of $N_2$.  Hence 
$\theta$ is an isomorphism.  
\end{proof}

Proposition \ref{hiProp1} implies that 
the isomorphism classes of $\rev_L$ form a {\em bona fide} 
partially ordered set.  There is a smallest element -
the class of the trivial lift $(\mathbb C, {\rm id}, \iota)$, $\iota: \mathbb C\to M$ denoting  
the inclusion $\iota(\lambda)=\lambda\cdot\mathbf 1_M$.  
The following result characterizes the class of an asymptotic 
lift as the largest element:  

\begin{prop}\label{hiThm1}
Let $(N_\infty,\alpha_\infty,E_\infty)$ be an asymptotic lift of $L$.  Then  
every $(N,\alpha,E)\in\rev_L$ satisfies $(N,\alpha,E)\leq (N_\infty, \alpha_\infty,E_\infty)$.  
\end{prop}

\begin{proof}  As in the proof of Theorem \ref{inThm1}, we can realize $N_\infty$ 
as the space of all diagonal operators $N_\infty=\{\dig(x_n): (x_n) \in\mathcal S\}$, 
$\alpha_\infty$ as the shift automorphism, and $E_\infty$ as the $0$th component map.  
Let $\theta: N\to N_\infty$ be the UCP map defined in (\ref{inEq0.1}).  
We have already pointed out in Remark \ref{inRem1} that $\theta$ is a 
homomorphism from $(N,\alpha,E)$ to $(N_\infty,\alpha_\infty,E_\infty)$.  
\end{proof}

\section{UCP maps on von Neumann algebras}\label{S:vn}

In this section we prove 
that the asymptotic lift of a UCP map acting on a von Neumann algebra is 
actually a W$^*$-dynamical system.  

\begin{thm}\label{vnThm1}
Let $(N,\alpha,E)$ be the asymptotic lift 
of a UCP map $L:M\to M$ that acts on a von Neumann algebra $M$.  Then 
$N$ is a von Neumann algebra and $\alpha$ is a $*$-automorphism of $N$.  
\end{thm}

We will deduce Theorem \ref{vnThm1} from the following proposition, 
which is normally used to 
establish the existence 
and uniqueness of the Poisson boundary of a noncommutative space 
of ``harmonic functions".  The noncommutative Poisson boundary is a far-reaching 
generalization of the fact that the space of bounded harmonic functions in 
the open unit disk $D$ is isometrically isomorphic to the abelian von Neumann 
algebra $L^\infty(\partial D, \frac{d\theta}{2\pi})$ 
of bounded measurable functions on the boundary $\partial D$ of $D$.  
We sketch a proof for the reader's convenience; 
more detail can be found in \cite{arvPoisBd} and \cite{izuPois}.  
The result itself appears to have been first discovered in \cite{EffSt}, Corollary 1.6.  

\begin{prop}\label{vnProp1}
Let $\Lambda: M\to M$ be a UCP map on a von Neumann algebra and let 
$H_\Lambda$ be the operator system of all harmonic elements of $M$ 
$$
H_\Lambda=\{x\in M: \Lambda(x)=x\}.  
$$
There is a unique associative multiplication $x,y\in H_\Lambda\mapsto x\circ y\in H_\Lambda$ that 
turns $H_\Lambda$ into a von Neumann algebra with predual $(H_\Lambda)_*$, 
on which the group $\aut H_\Lambda$ of automorphisms of the operator 
system structure of $H_\Lambda$ acts naturally as the group of all $*$-automorphisms.  
\end{prop}

\begin{proof}[Sketch of Proof] 
Uniqueness: Given two such multiplications, 
the identity map defines a complete order isomorphism of one $C^*$-algebra 
structure to the other.  Hence it is a $*$-isomorphism, and the 
two multiplications agree.  

Existence: 
We claim that there is a completely positive idempotent linear map $E:M\to M$ that has range 
$H_\Lambda$.  Indeed, if one topologizes the set of bounded linear maps 
of $M$ to itself with the topology of point-weak$^*$-convergence, then the set 
of completely positive unital maps on $M$ becomes a compact space.  Let 
$E$ be any 
limit point of the sequence of averages 
$$
A_n=\frac{1}{n}(\Lambda+\Lambda^2+\cdots+\Lambda^n),\qquad n=1,2,\dots.   
$$  
Using $\Lambda A_n=A_n\Lambda=\frac{n+1}{n}A_{n+1}-\frac{1}{n}\Lambda$, together with 
the straightforward estimate $\|A_n-A_{n+1}\|\leq \frac{2}{n+1}$, 
one finds that $\Lambda E=E\Lambda=E$, and from that follows $E^2=E$ as well as $E(M)=H_\Lambda$.  

Such an idempotent $E$ allows one to introduce a Choi-Effros multiplication on $H_\Lambda$, 
$x\circ y=E(xy)$, $x,y\in H_\Lambda$,  which makes $H_\Lambda$ into a $C^*$-algebra 
(\cite{ChE1}, Theorem 3.1, Corollary 3.2).  Since $H_\Lambda$ is weak$^*$-closed, it is 
the dual of the Banach space $(H_\Lambda)_*$, and a theorem of Sakai (\cite{sakaiC*W*}, 
Theorem 1.16.7) implies that
$H_\Lambda$ is a von Neumann algebra with predual $(H_{\Lambda})_*$.  
\end{proof}

\begin{proof}[Proof of Theorem \ref{vnThm1}] Let $(N,\alpha,E)$ be the 
asymptotic lift constructed in the proof of Theorem \ref{inThm1}, 
in which $N=\{\dig (x_n): (x_n)\in\mathcal S\}$, $\alpha$ is the bilateral shift 
automorphism $\alpha(X)=UXU^{-1}$, and $E(X)=X_{00}$, for $X\in\mathcal B(\tilde H)$, 
$\tilde H=\ell^2(\mathbb Z)\otimes H$ being the Hilbert space of sequences introduced 
there.  In order to prove the existence of a von Neumann algebra structure on $N$, we appeal to 
Proposition \ref{vnProp1} as follows.  
  
Consider the larger von Neumann algebra $\tilde M\supseteq N$  
of all bounded diagonal operators $Y=\diag (y_n)$ with 
components $y_n\in M$, and let $\Lambda$ be the map defined on $\tilde M$ by 
$$
\Lambda(\diag(y_n))=\diag(z_n), \qquad {\rm where\ } z_n=L(y_{n+1}),\quad n\in\mathbb Z.  
$$
Obviously, $\Lambda$ is a UCP map 
of $\tilde M$ with the property $\Lambda (Y)=Y$ if and only if $Y$ has the form 
$Y=\diag(x_n)$ with 
$(x_n)$ an inverse sequence for $L$.  Thus, $N$ is the space of all $\Lambda$-fixed 
elements of $\tilde M$.  An application of Proposition \ref{vnProp1} completes the 
proof.  
\end{proof}

\section{Nontriviality of the Asymptotic Dynamics}\label{S:nt}

A $W^*$-dynamical system $(A,\mathbb Z)$ is considered trivial if the 
automorphism of $A$ that implements the $\mathbb Z$-action is the identity automorphism.  
The purpose of this section is to show that the asymptotic lift of a UCP map is frequently  
a nontrivial dynamical system:   

\begin{thm}\label{ntThm1}
For every UCP map $L:M\to M$ of a dual operator system, the following are equivalent.
\begin{enumerate}
\item[(i)] The asymptotic lift $(N,\alpha,E)$ of $L$ satisfies $\alpha(y)=y$, $y\in N$.  
\item[(ii)]  Every operator $x$ in $\cap_nL^n(\ball M)$ satisfies $L(x)=x$.  

\item[(iii)]  The semigroup $\{L^n:n\geq 0\}$ oscillates slowly in the sense that  
$$
\lim_{n\to\infty}\|\rho\circ L^n-\rho\circ L^{n+1}\|=0,\qquad \rho\in M_*.  
$$
\end{enumerate}
\end{thm}

\begin{proof}
(i)$\implies$(ii):  We have $E(\ball N)=\cap_nL^n(\ball M)$ by Lemma \ref{inLem1}, so every 
$x\in\cap_nL^n(\ball M)$ has the form $x=E(y)$ for some $y\in N$.  Hence 
$L(x)=L(E(y))=E(\alpha(y))=E(y)=x$.  

(ii)$\implies$(iii): Fix $\rho\in M_*$.  Another application of Lemma \ref{inLem1} implies 
that $E(\ball N)=\cap_nL^n(\ball M)$ is pointwise fixed by $L$, 
hence $\|(\rho-\rho\circ L)\circ E\|=0$.  
Using 
(\ref{inEqA}), we obtain 
$$
\lim_{n\to\infty}\|\rho\circ L^n-\rho\circ L^{n+1}\|=\lim_{n\to\infty}\|(\rho-\rho\circ L)\circ L^n\|
=\|(\rho-\rho\circ L)\circ E\|, 
$$
and the right side is zero.  

(iii)$\implies$(i):  The preceding formula implies that $\|(\rho-\rho\circ L)\circ E\|=0$ 
for every $\rho\in M_*$, hence $E=L\circ E$.  So for every $n\geq 1$ and every $y\in N$, 
\begin{align*}
E(\alpha^{-n}(y-\alpha(y))&=E(\alpha^{-n}(y))-E(\alpha^{-n+1}(y))
\\&=E(\alpha^{-n}(y))-L(E(\alpha^{-n}(y)))=0
\end{align*}
and $\alpha(y)=y$ follows from nondegeneracy (\ref{inEq0}).  
\end{proof}

\begin{rem}[Matrix Algebras]  
Many UCP maps acting on matrix algebras 
are associated with nontrivial $W^*$-dynamical systems, because of the following observation:   
a UCP map $L$ of $M_n=M_n(\mathbb C)$ 
satisfies property (iii) of Theorem \ref{ntThm1}  $\iff$ $\sigma(L)\cap \mathbb T=\{1\}$.  
Indeed, to prove $\implies$ contrapositively, 
suppose there is a point $\lambda$ in the spectrum of $L$ 
that satisfies $\lambda\neq 1=|\lambda|$.  Choose 
an operator $x\in M_n$ satisfying $\|x\|=1$ and $L(x)=\lambda x$, and choose $\rho\in M_*$ 
such that $\rho(x)=1$.  Then for all $n\geq 0$ we have 
\begin{align*}
\|\rho\circ L^{n+1}-\rho\circ L^n\|&\geq |\rho(L^{n+1}(x))-\rho(L^n(x))|
=
|\lambda^{n+1}\rho(x)-\lambda^n\rho(x)|
\\&
=|\lambda-1|.   
\end{align*}
Hence 
$\{L^n\}$ does not oscillate slowly. 

Conversely, if $\sigma(L)\cap\mathbb T=\{1\}$, then since points of $\sigma(L)\cap\mathbb T$ 
are associated with {\em simple} eigenvectors of $L$, the spectrum of 
the restriction $L_0$ of $L$ 
to the range of ${\rm id}-L$ is contained in the open unit disk $\{z: |z|<1\}$, 
hence $\|L_0^n\|\to 0$ as $n\to\infty$.  It follows that $L^{n+1}-L^n$ tends 
to zero (in norm, say) as $n\to\infty$; and that obviously 
implies condition (iii) of Theorem \ref{ntThm1}.   We remark that the asymptotic behavior 
of $\|T^{n+1}-T^n\|$ for contractions $T$ on Banach spaces has been much-studied; 
see \cite{katzTz} and references therein.

Elementary examples show that 
any finite subset of the unit circle that contains $1$ and is stable under 
complex conjugation can occur as $\sigma(L)\cap\mathbb T$ for 
a UCP map $L$ of a matrix algebra $M_n$ (for appropriately large $n$).  
Hence there are many examples of 
UCP maps on finite-dimensional noncommutative von Neumann 
algebras whose asymptotic liftings are nontrivial finite-dimensional 
$W^*$-dynamical systems.  

The 
asymptotic behavior of UCP maps on finite-dimensional algebras is discussed more completely 
in Section \ref{S:ex}.  
\end{rem}

\begin{rem}[Automorphisms and Endomorphisms]
Automorphisms are at the opposite extreme from slowly oscillating 
UCP maps.  Indeed, if $\alpha$ is any automorphism of a dual operator system 
$M$, then $\alpha$ 
induces an isometry of the predual $M_*$ via $\omega\mapsto \omega\circ\alpha$, 
and for every $\rho\in M_*$ and $n\in\mathbb Z$ we have 
$$
\|\rho\circ\alpha^{n+1}-\rho\circ\alpha^n\|=\|(\rho\circ\alpha-\rho)\circ\alpha^n\|
=\|\rho\circ\alpha-\rho\|.  
$$
This formula obviously implies that the only slowly oscillating automorphism 
is the 
identity automorphism.  

If $\alpha$ is merely an isometric UCP map on $M$ and $M_\infty=\cap_n\alpha^n(M)$ is the 
``tail" operator system, then for every $\rho\in M_*$ we have 
$$
\|\rho\circ\alpha^{n+1}-\rho\circ\alpha^n\|=\|(\rho\circ\alpha-\rho)\circ\alpha^n\|=
\|(\rho\circ\alpha-\rho)\restriction_{\alpha^n(M)}\|, 
$$
and the right side decreases to 
$\|(\rho\circ\alpha-\rho)\restriction_{M_\infty}\|$ as $n\to\infty$ by Lemma \ref{inLem0}.  
Hence $\alpha$ oscillates slowly iff it  
restricts to the identity map on $M_\infty$.  
\end{rem}

\section{Identification of the asymptotic dynamics}\label{S:id}

Let $L:M\to M$ be a UCP map acting on a von Neumann algebra $M$ and let 
$(N,\alpha,E)$ be its asymptotic lift.  We have seen that $(N,\mathbb Z)$ is a $W^*$-dynamical 
system; indeed, the proof 
of Theorem \ref{vnThm1} shows that the algebraic structure of $N$ is that of 
the noncommutative Poisson boundary 
of an associated UCP map $\Lambda: \tilde M\to \tilde M$ on another von Neumann algebra.  

In general, it can be a significant problem to find a concrete realization of the 
Poisson boundary, even when $M=L^\infty(X,\mu)$ is commutative (see \cite{kaiVerRW} --
this is called the {\em identification problem} in \cite{kai1}).   
In the noncommutative case, there are only a few examples for which this 
problem has been effectively solved.  Three of them are discussed in \cite{izuPois}.  

In this section we contribute to this circle of ideas by showing that the 
asymptotic lift of a UCP map $L$ on a von Neumann algebra is isomorphic 
to the tail flow of the minimal dilation of $L$ to a $*$-endomorphism.  
We will clarify the precise relation 
between the asymptotic lift of $L$ and the Poisson boundary of $L$ in Section \ref{S:pb}.

It will not be necessary to reiterate   
explicit details of the dilation theory of UCP maps acting 
on von Neumann algebras  (see Chapter 8 of \cite{arvMono}).  
Instead, we simply recall that every 
UCP map $L:M\to M$ acting on a von Neumann algebra  
can be dilated {\em minimally} to a normal unit-preserving 
isometric $*$-endomorphism  of a larger von Neumann 
algebra $\alpha: N\to N$ that contains $M=pNp$ as a corner. 
Any two minimal dilations 
of $L$ are 
naturally isomorphic.  

More generally, any unit-preserving normal isometric $*$-endomorphism 
$\alpha: N\to N$ of a von Neumann algebra $N$ gives rise to a decreasing 
sequence of von Neumann subalgebras $N\supseteq \alpha(N)\supseteq \alpha^2(N)\supseteq\cdots$ 
whose intersection 
$$
A=\bigcap_{n\geq 0}\alpha^n(N)
$$
is a von Neumann algebra with the property that the restriction of $\alpha$ to $A$ is 
a {\em $*$-automorphism} of $A$.  That $W^*$-dynamical system 
$(A,\mathbb Z)$ (also written $(A,\alpha\restriction_A)$) 
is called the {\em tail flow} of the endomorphism $\alpha:N\to N$.  The 
tail flow  is clearly an interesting conjugacy 
invariant of the endomorphism 
$\alpha$, but it has received little attention in the past.

The following result clarifies the role of the tail flow in noncommutative 
dynamics by identifying 
the asymptotic lift of $L$ as the tail flow of the minimal dilation of $L$.  
Actually, we prove somewhat more, since the setting of Theorem \ref{isThm1} 
includes dilations that are not necessarily minimal.

\begin{thm}\label{isThm1}
Let $\alpha:N\to N$ be a unital normal isometric endomorphism of a von Neumann 
algebra $N$,  let $p\in N$ be a projection in $N$ satisfying $p\leq\alpha(p)$ 
and $\alpha^n(p)\uparrow \mathbf 1_N$ as $n\uparrow\infty$.  Let 
$M=pNp$ be the corresponding corner of $N$ and let $L:M\to M$ be the UCP map 
defined by 
$$
L^n(x)=p\alpha^n(x)p,\qquad x\in M,\quad n=0,1,2,\dots.     
$$

Then the asymptotic lift of $L:M\to M$ is isomorphic to $(A, \alpha\restriction_A,E)$, 
where $(A,\alpha\restriction_A)$ is the tail flow of $\alpha$ and 
$E: A\to M$ is the map $E(a)=pap$.
\end{thm}

\begin{proof} Obviously, $\alpha\restriction_A$ is an automorphism of the operator 
system structure of $A$, and $E$ is a UCP map of $A$ to $M=pNp$.  We claim that 
the nondegeneracy requirement $E(\alpha^{-n}(a))=0,\ n\geq 1\ \implies a=0$ is satisfied.
Indeed, fixing such an $a\in A$ and $n\geq 1$, we have 
$$
\alpha^n(p)a\alpha^n(p)=\alpha^n(p\alpha^{-n}(a)p)=\alpha^n(E(\alpha^{-n}(a))=0,\qquad n=1,2,\dots.  
$$
Since $\alpha^n(p)$ converges strongly to $\mathbf 1$ as $n\to \infty$, 
$\alpha^n(p)a\alpha^n(p)$ converges strongly to $a$, hence $a=0$.  

It remains to establish the formulas (\ref{inEqA}):
\begin{equation}\label{idEq1}
\lim_{k\to\infty}\|\rho\circ ({\rm id}_n\otimes L^k)\| = \|\rho\circ ({\rm id}_n\otimes E)\|, 
\end{equation}
for every $\rho\in(M_n\otimes M)_*$ and every $n=1,2,\dots$.  We 
first consider the case 
$n=1$.  Choose $\rho\in M_*$ and define $\bar\rho\in N_*$ by 
$\bar\rho(y)=\rho(pyp)$, $y\in N$.  For every $k\geq 1$, we claim
\begin{equation}\label{idEq2}
\|\rho\circ L^k\|=\|\bar\rho\circ\alpha^k\|,\quad {\rm and\ \quad}\|\rho\circ E\|=\|\bar\rho\restriction_A\|.  
\end{equation}
Indeed, since $L^k(x)=p\alpha^k(x)p$ for all $x\in M$, it follows that 
for every $y\in N$ we have 
$\rho(L^k(pyp))=\rho(p\alpha^k(pyp)p)=\rho(p\alpha^k(y)p)=\bar\rho(\alpha^k(y))$, 
and that identity clearly implies $\|\rho\circ L^k\|=\|\bar\rho\circ\alpha^k\|$.  The second 
formula of (\ref{idEq2}) follows from the identity 
$\rho(E(a))=\rho(pap)=\bar\rho(a)$ for $a\in A$.  

If we now apply (\ref{inEq6}) to the decreasing sequence of weak$^*$-compact sets 
$K_j=\alpha^j(\ball N)$, $j=1,2,\dots$, having intersection  $K_\infty=\ball A$, 
we find  that $\|\bar\rho\circ\alpha^j\|=\sup\{|\rho(y)|:y\in K_j\}$ and   
$\|\bar\rho\restriction_A\|=\sup\{|\bar\rho(a)|:a\in\ball A\}$, 
so by Lemma \ref{inLem0} we conclude  
$$
\lim_{k\to\infty}\|\rho\circ L^k\|=\lim_{k\to\infty}\|\bar\rho\circ\alpha^k\|=
\|\bar\rho\restriction_A\|=\|\rho\circ E\|, 
$$
proving (\ref{idEq1}).  

This argument applies {\em verbatim} to establish (\ref{idEq1})  
throughout the matrix hierarchy.  Indeed, for $n\geq 2$, the hypotheses of Theorem 
\ref{isThm1} carry over to the corner $M^{(n)}=p_nN^{(n)}p_n$, where $p_n=\mathbf 1_{M_n}\otimes p$, 
with $\alpha$ replaced by ${\rm id}_n\otimes \alpha$.  
We have $p_n\leq ({\rm id}_n\otimes \alpha)(p_n)\leq ({\rm id}_n\otimes \alpha^2)(p_n)\leq 
\cdots\uparrow \mathbf 1_{N^{(n)}}$, and for a fixed $\rho\in M^{(n)}_*$ there are 
appropriate versions of the formulas (\ref{idEq2}).  
\end{proof}

\section{Identification of the Poisson boundary}\label{S:pb}

In this section we show how the Poisson boundary can be 
characterized in terms of the asymptotic lift.  Since the asymptotic lift can often 
be calculated explicitly -- either directly as in Section \ref{S:ex} 
or using the tools of dilation theory and Theorem 
\ref{isThm1} -- this result contributes 
to the identification problem for noncommutative Poisson boundaries as 
discussed in 
Section \ref{S:id}.  

Let $L:M\to M$ be a UCP map acting on a von Neumann algebra, and consider the dual 
operator system of all noncommutative harmonic elements 
$$
H_L=\{x\in M: L(x)=x\}.  
$$
We have pointed out in Section \ref{S:vn} that $H_L$ carries a unique von Neumann 
algebra structure, and that von Neumann algebra is called the Poisson boundary of
the map $L:M\to M$.  Given a 
$*$-automorphism of a von Neumann algebra $N$, we write $N^\alpha=\{y\in N: \alpha(y)=y\}$ 
for its fixed subalgebra.

\begin{prop}
Let $L:M\to M$ be a UCP map on a von Neumann algebra, let $H_L$ be its 
Poisson boundary,  and let $(N,\alpha,E)$ be its 
asymptotic lift.  Then the restriction of $E$ to the fixed algebra $N^\alpha$ 
implements an isomorphism of von Neumann algebras $N^\alpha\cong H_L$.  
\end{prop}

\begin{proof}
The equivariance property $E\circ\alpha=L\circ E$, implies 
$E(N^\alpha)\subseteq H_L$.  
For the opposite inclusion, every element 
$a\in H_L$ gives rise to an inverse sequence $\bar a=(\cdots,a,a,a,\cdots)$ and, after realizing 
$(N,\alpha,E)$ concretely as in section \ref{S:al}, we conclude that 
$a=E(\dig \bar a)\in E(N^\alpha)$.  
Finally, a straightforward application of 
formula  (\ref{inEq5})  of Lemma \ref{inLem2} throughout the matrix hierarchy
shows that 
$E\restriction_{N^\alpha}$ 
is a complete isometry.  
Since both $N^\alpha$ and $H_L$ are von Neumann 
algebras, $E\restriction_{H_L}$ is a $*$-isomorphism.
\end{proof}

\section{Examples and Concluding Remarks}\label{S:ex}

Theorem \ref{isThm1} identifies the asymptotic lift of a UCP map in terms of 
its minimal dilation.  While one can often calculate properties of 
the minimal dilation in explicit terms 
(see Chapter 8 of \cite{arvMono}), those computations can be 
cumbersome and sometimes difficult.  
On the other hand, given a specific 
UCP map, we have 
found that it is 
often easier to calculate 
its asymptotic lift directly in concrete terms.  The purpose of this section 
is to illustrate that fact by carrying out calculations for some 
examples that require a variety of techniques.  

\subsection{Stochastic Matrices}\label{SS:st}
It is appropriate to begin with the classical commutative case 
having its origins in the theory of Markov chains.  
Let $P=(p_{ij})$ be an $n\times n$ matrix of nonnegative numbers satisfying 
$\sum_j p_{ij}=1$ for every $i=1,\dots,n$.  If we view the elements of 
the von Neumann algebra $M=\mathbb C^n$ as column vectors, then $P$ gives 
rise to a UCP map on $M$ by matrix multiplication.  We now calculate the 
asymptotic lift of $P$, and we relate that to classical results of Frobenius 
\cite{frob1}, \cite{frob2}, generalizing earlier results of Perron \cite{perron1}, \cite{perron2}, 
on the asymptotic behavior of such matrices.  To keep the discussion as simple 
as possible, we restrict 
attention to the case where $P$ is {\em irreducible} in the sense that 
the only projections $e\in M$ satisfying $P(e)\leq e$ are $e=0$ and $e=\mathbf 1$ 
(but see Remark \ref{stRem1}).  
We write $\sigma(P)$ for the spectrum of the linear 
operator $P:\mathbb C^n\to \mathbb C^n$;  $\sigma(P)$ is
a subset of the closed unit disk that contains the eigenvalue $1$.  

In this context, Theorem 2 of (\cite{gant}, see pp. 65--75) can be paraphrased 
as follows.
\begin{thm}\label{exThm2}
Let $P$ be an irreducible stochastic $n\times n$ matrix.   Then there is a $k$,
$1\leq k\leq n$ such that $\sigma(P)\cap\mathbb T=\{1, \zeta,\zeta^2,\dots,\zeta^{k-1}\}$ is 
the set of all $k$th roots of unity, 
$\zeta=e^{2\pi i/k}$, each $\zeta^j$ being a simple eigenvalue.  When 
$k>1$, there is 
a permutation matrix $U$ such that $UPU^{-1}$ has cyclic form 
\begin{equation}\label{stEq1}
UPU^{-1}=
\begin{pmatrix}
0&C_0&0&\dots&0\\
0&0&C_1&\dots&0\\
\vdots&\vdots&\vdots&\dots&\vdots\\
0&0&0&\dots&C_{k-2}\\
C_{k-1}&0&0&\dots&0
\end{pmatrix}
\end{equation}
in which the $C_j$ are rectangular submatrices, 
and where the zero submatrices along the diagonal are all square.  
\end{thm}

Notice that we can eliminate the unitary permutation matrix $U$ entirely by 
replacing $P$ with a suitably relabeled stochastic matrix, and we do so.  
For each $j=0,1,\dots,k-1$, let $e_j$ be the projection corresponding 
to the domain subspace of the block $C_j$.  These projections are mutually 
orthogonal, have sum $\mathbf 1$, and satisfy $P(e_j)\leq e_{j+1}$ (addition modulo $k$).  
Making use of $P(\mathbf 1)=\mathbf 1$, we find that in fact, $P(e_j)=e_{j+1}$ for every $j$.  
Let $N$ be the linear span of $e_0,\dots,e_{k-1}$ and let $\alpha$ be the 
restriction of $P$ to $N$.  Clearly $N$ is a $*$-subalgebra of $M$ 
and $\alpha:N\to N$ is the $*$-automorphism associated with this 
cyclic permutation of the minimal projections $e_0,e_1,\dots,e_{k-1}$ of $N$.  

\begin{prop}\label{exProp1}
The asymptotic lift of $P:M\to M$ is the triple $(N,\alpha,E)$, 
where 
$E:N\subseteq M$ is the inclusion map.  
\end{prop}

\begin{proof}  For $k>1$ as above, consider the UCP map 
$P^k: M\to M$.  The spectrum of $P^k$ consists of the 
eigenvalue $1$, together with other spectral points in the open unit 
disk $\{|z|<1\}$, the eigenvalue $1$ having eigenspace $N$.  Hence the sequence 
of powers $P^k, P^{2k}, P^{3k},\dots$ converges to an idempotent  
linear map $Q$ on $M$ having range $N$.  Clearly $Q$ is a UCP projection 
onto $N$.  For every $\rho\in M_*$, the norms $\|\rho\circ P^r\|$ decrease 
as $r$ increases, hence 
\begin{equation}\label{stEq2}
\lim_{r\to\infty}\|\rho\circ P^r\|=\lim_{m\to\infty}\|\rho\circ (P^{k})^m\|=
\|\rho\circ Q\|=\|\rho\restriction_N\|=\|\rho\circ E\|.  
\end{equation}
The same argument applies throughout the matrix hierarchy over $M$, and 
we conclude that the two criteria of Definition \ref{inDef1} are satisfied.  
Hence $(N,\alpha,E)$ is the asymptotic lift of $P$.  
\end{proof}

\begin{rem}[Asymptotics of the powers of $P$]\label{stRem2}  The idempotent UCP map 
$$
Q=\lim_{m\to\infty}P^{mk}
$$ 
exhibited in the proof of Proposition \ref{exProp1} provides 
a precise sense in which the asymptotic lift $(N,\alpha,E)$ 
contains all of the asymptotic information about the sequence of 
powers $P, P^2, P^3,\dots$.
Indeed, we   
claim that there are positive constants $c, r$ with $ r<1$ such that 
\begin{equation}\label{stEq3}
\|P^n-\alpha^nQ\|\leq cr^n,\qquad n=1,2,3,\dots.  
\end{equation}
This follows from the fact that $Q$ is an idempotent that commutes with 
$P$ with the property that $\alpha=P\restriction_{N=\ran Q}$ has spectrum in the unit circle 
and $P\restriction_{\ran ({\rm id}-Q)}$ has spectrum in $\{|z|<1\}$.  
Choosing a number $r<1$ so that $\{|z|<r\}$ contains the spectrum of 
$P\restriction_{\ran ({\rm id}-Q)}$, the spectral radius formula of elementary 
Banach algebra theory implies that the sequence 
$$
r^{-n}\|P^n-\alpha^nQ\|=r^{-n}\|P^n-P^nQ\|=(r^{-1}\|P^n\restriction_{\ran({\rm id}-Q)}\|^{1/n})^n
$$
tends to zero as $n\to \infty$, hence there is a $c>0$ so that (\ref{stEq3}) 
is satisfied.  
\end{rem}

\begin{rem}[Reducibility and Noncommutativity]\label{stRem1}
In the following subsection, we generalize the above 
result to the case of UCP maps acting on noncommutative finite-dimensional 
von Neumann algebras.  In particular, the discussion of Subsection \ref{SS:fd} 
applies equally 
to the reducible cases not covered in the statement of Theorem \ref{exThm2}.  
\end{rem}

\subsection{Finite-dimensional von Neumann algebras}\label{SS:fd}
One can compute the asymptotic lift of a UCP map 
on a finite-dimensional von Neumann algebra explicitly, using nothing but elementary 
methods along with the Choi-Effros multiplication.  Since these results 
have significance for quantum computing \cite{kup1}, we sketch that calculation  
in some detail.  

Let $L:M\to M$ be 
a UCP map on a finite-dimensional von Neumann algebra.  For every point  $\lambda\in\sigma(L)\cap \mathbb T$ let 
$$
N_\lambda=\{x\in M: L(x)=\lambda x\}  
$$
and let 
$$
N=\sum_{\lambda\in\sigma(L)\cap\mathbb T} N_\lambda
$$
be the sum of these maximal eigenspaces.  The identity operator belongs to $N$, 
and from the property $L(x^*)=L(x)^*$, $x\in M$, one deduces $N^*=N$.  Hence $N$ 
is an operator system such that the restriction $\alpha=L\restriction_N$ 
of $L$ to $N$ is 
a {\em diagonalizable} UCP map with spectrum $\sigma(L)\cap\mathbb T$.  Indeed, 
it is not hard to show that $\alpha$ is a UCP automorphism of $N$.

We digress momentarily to point out that in the classical setting of 
Frobenius' result for irreducible stochastic matrices $P$ as formulated in 
Theorem \ref{exThm2}, $N$ coincides with the span of the projections $e_0,\dots,e_{k-1}$ 
constructed above, the eigenvector associated with a 
$k$th root of unity $\lambda\in \sigma(P)\cap\mathbb T$ 
being given by 
$$
x_\lambda=e_0+\bar\lambda e_1+\bar\lambda^2 e_2+\cdots+\bar\lambda^{k-1}e_{k-1}.  
$$
It follows that, in such commutative cases, $N$ is closed under multiplication.

In the more general setting under discussion here, $\sigma(L)\cap \mathbb T$ need 
not consist of roots of unity and 
$N$ need not be closed under multiplication.  But in all cases $N$ 
can be made into a von Neumann 
algebra.  The most transparent proof of that fact uses
the following observation of Kuperberg 
(\cite{kup1}, also see Theorem 2.6 of \cite{arvStableI}):

\begin{lem}\label{fdLem1}
There is 
an increasing sequence of integers $n_1<n_2<\cdots$ 
such that $L^{n_1}, L^{n_2},\dots$ converges to an idempotent UCP map $Q$ 
with the property $N=Q(M)$.  In fact, $Q$  
is the unique idempotent limit point of the sequence of powers $L,L^2,L^3,\dots$.  
\end{lem}

One now uses the idempotent $Q$ to introduce a Choi-Effros multiplication 
$$
x\circ y=Q(xy), \qquad x,y\in N,
$$ 
in $N$, thereby making it into  a finite-dimensional 
von Neumann algebra.  Hence $(N,\alpha)$ becomes a (typically noncommutative) 
W$^*$-dynamical system, and the inclusion $E:N\subseteq M$ becomes a UCP map 
satisfying $E\circ\alpha=L\circ E$.

\begin{prop}
The triple $(N,\alpha,E)$ is the asymptotic lift of $L:M\to M$.  
\end{prop}

\begin{proof}[Sketch of proof]
The inclusion of $N$ in $M$ is injective, hence the nondegeneracy condition 
(\ref{inEq0}) is satisfied.  Thus we need only show that equality holds in the 
formulas (\ref{inEqA}), and that follows by a proof paralleling that of (\ref{stEq2}).
Indeed, letting $n_1<n_2<\cdots$ be a sequence such that $L^{n_k}\to Q$ 
as in Lemma \ref{fdLem1},   one finds 
that for every bounded linear functional $\rho$ on $M$, 
$$
\lim_{n\to\infty}\|\rho\circ L^n\|=\lim_{k\to\infty}\|\rho\circ L^{n_k}\|=
\|\rho\circ Q\|=\|\rho\restriction_N\|=\|\rho\circ E\|.  
$$
Obviously, the same argument can be promoted throughout the matrix hierarchy, 
as one allows $\rho$ to range  
over the dual of $M^{(n)}$, $n=1,2,\dots$.  
\end{proof}

\begin{rem}[Asymptotics of the powers of $L$.]
Making use of the idempotent $Q$ much as in Remark \ref{stRem2}, 
one finds that $(N,\alpha,E)$ contains 
all asymptotic information about the sequence $L, L^2, L^3,\dots$ 
because of the following precise 
estimate: There are positive constants $c, r$ such that $r<1$ and 
\begin{equation}\label{fdEq1}
\|L^n-\alpha^nQ\|\leq cr^n,\qquad n=1,2, \dots.  
\end{equation}
\end{rem}

\subsection{UCP maps on $II_1$ factors.}\label{SS:er}

We now  calculate the asymptotic lifts of a family of nontrivial 
UCP maps acting on the hyperfinite $II_1$ factor $R$, the point being to show that 
the asymptotic lifts of these maps 
can be arbitrary $*$-automorphisms of $R$.  There are many variations 
on these examples that exhibit a variety of phenomena in 
other von Neumann algebras.  Here, 
we confine ourselves to the simplest nontrivial cases.  

Let $\tau$ be the tracial state of $R$, and let $P$ be any normal UCP map of 
$R$ having $\tau$ as an absorbing state in the sense that for every normal state 
$\rho$ of $R$, one has 
\begin{equation}\label{erEq1}
\lim_{k\to\infty}\|\rho\circ P^k-\tau\|=0.  
\end{equation}
Of course, the simplest such map is $P(x)=\tau(x)\mathbf 1$; but we have less trivial 
examples in mind.  For example, let $A_\theta=C^*(U,V)$ be the irrational rotation 
$C^*$-algebra, where $U,V$ are unitaries satisfying 
$UV=e^{2\pi i\theta}VU$ and $\theta$ is an irrational number.  
Then for every number 
$\lambda$ in the open unit interval $(0,1)$ one can 
show that there is a completely positive unit-preserving map $P_\lambda$ 
that is defined 
uniquely on $A_\theta$ by its action on monomials: 
$$
P_\lambda(U^pV^q)=\lambda^{|p|+|q|}U^pV^q,\qquad p,q\in\mathbb Z.   
$$  
Moreover, since in every finite linear combination $x=\sum_{p,q}a_{p,q}U^pV^q$, 
the coefficients satisfy 
$|a_{p,q}|\leq \|x\|$, $p,q\in\mathbb Z$, a straightforward estimate leads to 
$$
\sup_{\|x\|\leq 1}\|P_\lambda^k(x)-\tau(x)\mathbf 1\|\leq\sum_{(p,q)\neq (0,0)}\lambda^{(|p|+|q|)^k}
\leq \lambda^{k-1}\sum_{(p,q)\neq(0,0)}\lambda^{|p|+|q|}=c\lambda^{k},   
$$
for some positive constant $c$.   We conclude that 
$$
\lim_{k\to\infty}\|P^k_\lambda-\tau(\cdot)\mathbf 1\|=0,   
$$ 
and in particular, (\ref{erEq1}) holds for 
every state $\rho$  of $C^*(U,V)$.
  
The GNS construction applied to the tracial state of $A_\theta$ 
now provides a representation $\pi$ of $A_\theta$ with the property that the weak 
closure of $\pi(A_\theta)$ is $R$, and there is a unique UCP map $P$ on 
$R$ such that $P(\pi(x))=\pi(P_\lambda(x))$ for all $x\in A_\theta$.  
The preceding paragraph implies that the extended 
map $P:R\to R$ has the asserted property (\ref{erEq1}).  

Similar examples of UCP maps acting on other $II_1$ factors can be constructed 
by replacing $R$ with the group von Neumann algebras of discrete groups with infinite 
conjugacy classes (see Proposition 4.4 of \cite{arvStableI}).  

Choose a UCP map $P:R\to R$ having property (\ref{erEq1}), 
choose an arbitrary $*$-automorphism $\alpha$ of $R$, 
and consider the UCP map $L=P\otimes \alpha$ defined uniquely on the spatial 
tensor product $R\otimes R$ by 
$$
L(x\otimes y)=P(x)\otimes \alpha(y),\qquad x,y\in R.  
$$
Since $R\otimes R$ is isomorphic to $R$, one can think of $L$ as a UCP map 
on $R$.  

\begin{prop}\label{erProp1}
Let $E: R\to R\otimes R$ be the map $E(x)=\mathbf 1\otimes x$, $x\in R$.  
The asymptotic lift of $L=P\otimes\alpha$ is the triple $(R,\alpha,E)$.  
\end{prop}

\begin{proof}[Sketch of proof]  $E$ is clearly an injective UCP map of von 
Neumann algebras satisfying the equivariance condition $E\circ\alpha=L\circ E$.  
Thus it remains only to verify the formulas (\ref{inEqA}).  For that, we will prove 
a stronger asymptotic relation.  
Let $Q: R\otimes R\to \mathbf 1\otimes R$ be the 
$\tau$-preserving conditional expectation 
$$
Q(x\otimes y)=\tau(x)y,\qquad x,y\in R.  
$$
We claim that for all $\rho\in (R\otimes R)_*$, one has 
\begin{equation}\label{erEq2}
\lim_{k\to\infty}\|\rho\circ L^k-\rho\circ ({\rm id_R}\otimes\alpha)^k\circ Q\|=0.
\end{equation}
Indeed, since $(R\otimes R)_*$ is the norm-closed linear span of functionals of the 
form $\rho_1\otimes \rho_2$ with $\rho_k\in R_*$, obvious estimates show that 
it suffices to prove 
(\ref{erEq2}) for decomposable functionals of the form $\rho_1\otimes\rho_2$, where 
$\rho_1,\rho_2\in R_*$.  Fix $\rho=\rho_1\otimes \rho_2$ of this form.  
Using the decomposition 
$$
(L^k-({\id_R}\otimes \alpha)^k\circ Q)(x\otimes y)=(P^k(x)-\tau(x)\mathbf 1)\otimes\alpha^k(y), 
$$
we find that $\|\rho_1\otimes\rho_2  (L^k-({\id_R}\otimes \alpha)^k\circ Q)\|$ decomposes 
into a product
\begin{align*}
\|(\rho_1\otimes\rho_2)((P^k-\tau(\cdot)\mathbf 1)\otimes\alpha^k)\|
&=
\|\rho_1\circ(P^k-\tau(\cdot)\mathbf 1)\|\cdot\|\rho_2\circ\alpha^k\|
\\&=
\|\rho_1\circ(P^k-\tau(\cdot)\mathbf 1)\|\cdot\|\rho_2\|
\\&=\|\rho_1\circ P^k-\tau\|\cdot\|\rho_2\|, 
\end{align*}
which by (\ref{erEq1}), tends to zero as $k\to\infty$. 

Note that (\ref{erEq2}) leads immediately to the case $n=1$ of (\ref{inEqA}), since 
$$
\lim_{k\to\infty}\|\rho\circ L^k\|=\lim_{k\to\infty}\|\rho\circ({\rm id}_R\otimes\alpha^k)\circ Q\|
=\lim_{k\to\infty}\|\rho\circ Q\circ\alpha^k\|=\|\rho\circ Q\|=\|\rho\circ E\|.  
$$

With trivial changes, these arguments can be repeated throughout the matrix 
hierarchy, after one identifies $M_n\otimes (R\otimes R)$ with $(M_n\otimes R)\otimes R$.  
We omit those mind-numbing details.  
\end{proof}

\subsection{The CCR heat flow}

The three asymptotic assertions (\ref{stEq3}), (\ref{fdEq1}), (\ref{erEq2}) 
are much 
stronger than the requirements of (\ref{inEqA}), and one might ask if those 
stronger results can be established for the asymptotic lifts $(N,\alpha,E)$ of more general UCP 
maps on von Neumann algebras $L:M\to M$.  
In each of the preceding examples, it was possible to identify $N$ with a dual operator 
subsystem $M_\infty\subseteq M$ (namely the range $E(N)$ of $E$), $\alpha$ with the restriction 
$\alpha_\infty=L\restriction_{M_\infty}$ of $L$ to $N$, and $E$ with the 
inclusion map $\iota: N\subseteq M$.  
There was also 
an idempotent completely positive projection $Q$ mapping $M$ onto $M_\infty$, and together,  these 
objects gave rise to the asymptotic relations 
\begin{equation}\label{lastEq}
\lim_{k\to\infty}\|\rho\circ L^k-\rho\circ \alpha_\infty^k\circ Q\|=0,\qquad \rho\in M_*.  
\end{equation}
The relative strength of (\ref{lastEq}) and (\ref{inEqA}) is clearly seen if one 
reformulates the case $n=1$ of (\ref{inEqA}) in this context as the following assertion
$$
\lim_{k\to\infty}|\,\|\rho\circ L^k\|-\|\rho\circ \alpha_\infty^k\circ Q\|\,|=
\lim_{k\to\infty}|\,\|\rho\circ L^k\|-\|\rho\circ E\|\,|=0, 
$$
(see the proof of Proposition \ref{erProp1}).  

So it is natural to ask if the stronger 
relations (\ref{lastEq}) can be established more generally, at least 
in cases where $N=M_\infty\subseteq M$ is a subspace of $M$ 
and $\alpha=\alpha_\infty$ is obtained by restricting $L$ to $M_\infty$.  
That is true, for example, in the more restricted context of \cite{arvStableI}.  
The purpose of these remarks is to show that the answer is no in general, by 
describing examples 
with the two 
properties $M_\infty\subseteq M$ and $\alpha=L\restriction_{M_\infty}$, but for which 
there is no idempotent completely positive map $Q$ satisfying (\ref{lastEq}).

The CCR heat flow is a semigroup of UCP maps $\{P_t:t\geq 0\}$ acting on 
the von Neumann algebra $M=\mathcal B(H)$
\cite{ccrHeat}.  Consider the single UCP map $L=P_{t_0}$ for some fixed 
$t_0>0$.  This map has the following two properties: a) there is no normal state 
$\rho\in M_*$ satisfying $\rho\circ L=\rho$, and b) for any two normal 
states $\rho_1, \rho_2\in M_*$ one has 
\begin{equation}\label{erEq6}
\lim_{k\to\infty}\|\rho_1\circ L^k-\rho_2\circ L^k\|=0.  
\end{equation}
We claim first that the asymptotic lift of this 
$L$ is the triple $(\mathbb C, {\rm id}, \iota)$, 
where $\iota$ is the inclusion of $\mathbb C$ in $M$, $\iota(\lambda)=\lambda\cdot\mathbf 1_M$.  
Indeed, to 
sketch the proof of the key assertion -- namely the case $n=1$ of formula (\ref{inEqA}) -- choose 
$\rho\in M_*$.   We have $\|\rho\circ \iota\|=|\rho(\mathbf 1)|$, so that (\ref{inEqA}) reduces in 
this case to the formula 
\begin{equation}\label{erEq7}
\lim_{k\to\infty}\|\rho\circ L^k\|=|\rho(\mathbf 1)|, \qquad \rho\in M_*.  
\end{equation}
It is a simple exercise to show that (\ref{erEq7}) and (\ref{erEq6}) 
are in fact equivalent assertions about $L$ -- note for example that (\ref{erEq6}) is the special 
case of (\ref{erEq7}) in which $\rho(\mathbf 1)=0$ -- and of course 
the restriction of $L$ to $\mathbb C\cdot \mathbf 1$ 
is the identity map $\alpha_\infty={\rm id}$.  This argument 
promotes naturally throughout the matrix hierarchy over $M$, hence 
$(\mathbb C,{\rm id},\iota)$ is the asymptotic lift of $L$.  

We now show that one cannot obtain formulas 
like (\ref{lastEq}) for this example.

\begin{prop}
There is no completely positive projection $Q$ of $M$ on $\mathbb C\cdot\mathbf 1$
that satisfies 
\begin{equation}\label{erEq8}
\lim_{k\to\infty}\|\rho\circ L^k-\rho\circ L^k\circ Q\|=0, \qquad \rho\in M_*.  
\end{equation}
\end{prop}

\begin{proof}
Indeed, a completely positive idempotent 
$Q$ with range $\mathbb C\cdot\mathbf 1$ 
would have the form $Q(x)=\omega(x)\mathbf 1$, $x\in M$,  
where $\omega$ is a state of $M$.  
For any normal state $\rho$ we have $\rho\circ L^k\circ Q(x)=
\rho(L^k(\omega(x)\mathbf 1))=\omega(x)$, 
$x\in M$, hence in this case (\ref{erEq8}) makes the assertion 
$$
\lim_{k\to\infty}\|\rho\circ L^k-\omega\|=0;    
$$
i.e., $\omega$ is an {\em absorbing state} for $L$.  But an absorbing state $\omega$ for 
$L$ is a {\em normal} state that satisfies $\omega\circ L=\omega$, contradicting 
property 
a) above.  
\end{proof} 

\begin{rem}[Further examples]
One can generalize this construction based on $L$.  For example, let $\alpha$  be a 
$*$-automorphism of $\mathcal B(H)$ and consider the the UCP map $L\otimes \alpha$ defined on 
$M=\mathcal B(H\otimes H)$.  One can show that $L\otimes \alpha$ 
has asymptotic lift $(\mathcal B(H),\alpha,E)$ 
where $E:\mathcal B(H)\to M$ is the UCP map $E(x)=\mathbf 1\otimes x$, much as in the 
proof of Proposition \ref{erProp1}.  With suitable choices of $\alpha$ (specifically, 
for any 
$\alpha$ that has a normal invariant state), one can 
also show that there is no completely positive projection $Q: M\to \mathbf 1\otimes\mathcal B(H)$
that satisfies (\ref{erEq8}).  Thus, even though a UCP map on a von Neumann 
algebra always has an asymptotic lift, {\em there are many examples for which one cannot 
expect precise asymptotic formulas such as (\ref{lastEq})}.  
\end{rem}

\bibliographystyle{alpha}

\newcommand{\noopsort}[1]{} \newcommand{\printfirst}[2]{#1}
  \newcommand{\singleletter}[1]{#1} \newcommand{\switchargs}[2]{#2#1}

\end{document}